\documentclass{article}
\usepackage{amsmath,amssymb,amsthm}
\DeclareMathOperator{\lcm}{lcm}
\DeclareMathOperator{\ord}{ord}
\newtheorem{theorem}{Theorem}
\newtheorem*{midy}{Midy's Theorem}
\newtheorem*{remark}{Remark}

\begin{document}


\title{\bf Midy's Theorem for Periodic Decimals}
\author{Joseph Lewittes}
\date{}


\maketitle


\section{Introduction}\label{s1}


It is well known--and a proof will appear in our subsequent discussion--that
 any rational number $c/d$, with $d$ relatively prime to $10$, has a purely
 periodic decimal expansion of the form
 $I.a_1 a_2 \ldots a_n a_1 \ldots a_n a_1 \ldots$, where $I$ is an integer,
 $a_1$, $a_2$, $\ldots$, $a_n$ are digits, and the block $a_1 a_2 \ldots a_n$
 repeats forever.
The repeating block is called the {\em period} and $n$ is its {\em length}.
We write the decimal as $I.\overline{a_1 a_2 \ldots a_n}$, the bar indicating
 the period.
Consider a few examples: $1/3 = 0.\overline{3}$, $1/7 = 0.\overline{142857}$,
 $2/11 = 0.\overline{18}$, $1/13 = 0.\overline{076923}$,
 $2/13 = 0.\overline{153846}$, $1/17 = 0.\overline{0588235294117647}$,
 $1/37 = 0.\overline{027}$, $1/73 = 0.\overline{01369863}$.
Note that when the period length is even and the period is broken into two
 halves of equal length which are then added, the result is a string of 9's.
Thus $142 + 857 = 999$, $1 + 8 = 9$, $076 + 923 = 999$, and so on; the
 numerator plays no role.
In each of these examples the denominator is a prime number.
Try a few composite denominators: $77 = 7 \times 11$,
 $1/77 = 0.\overline{012987}$; $803 = 11 \times 73$,
 $1/803 = 0.\overline{00124533}$; $121 = 11 \times 11$,
 $1/121 = 0.\overline{0082644628099173553719}$.
We see the property holds for $77$ and $121$ but fails for $803$.
According to Dickson \cite[p. 161, footnote 19]{1}, H. Goodwyn was apparently
 the first to observe (in print, 1802) this phenomenon for prime denominators,
 based on experimental evidence.
Over the past two centuries it has been rediscovered many times; it is called
 the `nines property' by Leavitt \cite{4} and `complementarity' by
 Shrader-Frechette \cite{7}.
This latter reference contains a historical perspective and a bibliography of
 the topic.


In 1836 E. Midy \cite{6} published at Nantes, France, a pamphlet of twenty-one
 pages on some topics in number theory with applications to decimals.
He was the first to actually prove something about our topic.
We formulate our own version of his main result.
As usual, $\gcd(a,b)$ denotes the greatest common divisor of the integers
 $a$, $b$.


\begin{midy}
Let $x$ and $N$ be positive integers, with $N > 1$, $\gcd(N,10) = 1$,
 $\gcd(x,N) = 1$ and $1 \le x < N$.
Assume $x/N = 0.\overline{a_1 a_2 \ldots a_{2k}}$ has even period length $2k$.
If
\begin{enumerate}
  \item[(i)] $N$ is a prime, or
  \item[(ii)] $N$ is a prime power, or
  \item[(iii)] $\gcd(N,10^k - 1) = 1$,
\end{enumerate}
then, for $1 \le i \le k$,
  \begin{equation}\label{eq1}
    a_i + a_{k + i} = 9.
  \end{equation}
\end{midy}


We refer to (\ref{eq1}) as the Midy property for the denominator $N$.
Proof of the theorem will be included later--see the Remark after Theorem 6 in
 Section 3--as part of a more general theory.
Note that $(iii)$ explains the difference between $1/77$ where $k = 3$ and
 $\gcd(77,10^3 - 1) = 1$, which has the Midy property, and $1/803$ where
 $k = 4$ and $\gcd(803,10^4 - 1) = 11$, for which the Midy property fails.


Various authors have given proofs of this theorem, or parts of it, most being
 unaware of Midy; even those who do cite him do so only through Dickson's
 reference \cite[p. 163]{1}.
Undoubtedly this is due to the obscure publication of Midy's paper.
Recently, Ginsberg \cite{2} extended Midy's theorem to the case where the
 period has length $3k$; he showed that when the period is broken into three
 pieces of length $k$ each and then added, the sum is again a string of $9$'s.
However, his result is stated only for fractions $1/p$, $p$ a prime, and
 numerator restricted to be $1$.
Example: $1/13 = 0.\overline{076923}$, $07 + 69 + 23 = 99$.
However, note that $2/13 = 0.\overline{153846}$, $15 + 38 + 46 = 99$,
 $3/13 = 0.\overline{230769}$, $23 + 07 + 69 = 99$,
 $1/21 = 0.\overline{047619}$, $04 + 76 + 19 = 99$, all of which suggest a
 wider application of the result.
This will be discussed in the final section.


Eventually, I decided to actually look at Midy's paper--it is available on
 microfilm at the New York Public Library--and, remarkably, Midy's approach
 enables one to prove a general theorem that includes the above results and
 even more.
Midy himself considered the case of period length $3k$, but he focused on the
 sums $a_i + a_{i + k} + a_{i + 2k}$, $1 \le i \le k$, which do not give smooth
 results.
For example, with $1/7$ as above, $3k = 6$, $k = 2$, $1 + 2 + 5 = 8$,
 $4 + 8 + 7 = 19$, even though $14 + 28 + 57 = 99$.
In fact, one easily sees that for period length $2k$ the two halves adding up
 to a string of $9$'s is equivalent to $a_i + a_{k + i} = 9$, $1 \le i \le k$,
 but for length $3k$ it is not so, as carrying may occur.


In the next section we concentrate on when the period can be broken up into
 blocks of equal length $k$ that add up to a multiple of $10^k - 1$ and, in the
 final section, we return to the question as to when the sum is exactly
 $10^k - 1$, a string of $k$ $9$'s.
Since it is just as easy to carry out the analysis for an arbitrary number base
 $B$ as for the decimal base $10$, we do so.


\section{Base B and Midy}\label{s2}


Let $B$ denote an integer $> 1$ which will be the base for our numerals.
The digits in base $B$, $B$-digits for short, are the numbers $0$, $1$, $2$,
 $\ldots$, $B - 1$.
Every positive integer $c$ has a unique representation as
 $c = d_{n - 1}B^{n - 1} + d_{n - 2}B^{n - 2} + \ldots + d_1B + d_0$, where $n$
 is a positive integer, each $d_i$ is a $B$-digit and $d_{n-1} > 0$.
As in the decimal case, where $B = 10$, we write $c$ in base $B$ as the numeral
 $d_{n-1} d_{n-2} \ldots d_0$.
When necessary to indicate the base, we write $[d_{n-1} d_{n-2} \ldots d_0]_B$.
For $B = 10$ we use the usual notation.
We now fix some notation.
Unless otherwise noted, our variables $a$, $b$, $\ldots$ denote positive
 integers.
$a|b$ indicates $a$ divides $b$. $B$ is the base and $N$, which will
 be the denominator of our fractions, is relatively prime to $B$.
$N^*$ is the set $\{x | 1 \le x \le N$ and $\gcd(x,N) = 1\}$, the set of
 positive integers less than $N$ relatively prime to $N$.
These will be the numerators of our fractions.
For $x \in N^*$, $x/N$ is a reduced fraction strictly between $0$ and $1$ and
 we are interested in the base $B$ expansion of such a fraction.
Recalling the elementary school long division process for the decimal expansion
 of fractions one sees that it amounts to the following.
Set $x_1 = x$, let $a_1$  be the integer quotient and $x_2$ the remainder when
 $Bx_1$ is divided by $N$.
Thus $Bx_1 = a_1N + x_2$, $0 \le x_2 < N$ and $a_1 = \lfloor Bx_1/N\rfloor$
 where $\lfloor~\rfloor$ is the greatest integer, or floor, function.
Continuing inductively, we obtain the following infinite sequence of equations,
 which we call the long division algorithm.


\begin{equation}\label{eq2}
  \begin{matrix}
    Bx_1 = a_1 N + x_2       \\
    Bx_2 = a_2 N + x_3       \\
    \dots                    \\
    Bx_i = a_i N + x_{i + 1} \\
    \dots
  \end{matrix}
\end{equation}


Since $0 < x_1/N < 1$, $Bx_1/N < B$, $a_1 = \lfloor Bx_1/N \rfloor < B$, so
 $a_1$ is a $B$-digit.
Also $B$ and $x_1$ are both relatively prime to $N$ so
 $Bx_1 \equiv x_2 \pmod N$ shows $(x_2,N) = 1$, so $x_2 \in N^*$.
In the same way, for all $i \ge 1$, $a_i$ is a $B$-digit and $x_i \in N^*$.
Dividing the first equation by $BN$, the second by $B^2N$, and in general the
 $i$th by $B^iN$ shows
 $x_1/N = a_1/B + a_2/B^2 + \ldots + a_i/B^i + x_{i + 1}/B^iN$.
Since $0 < x_{i + 1}/B^iN < 1/B^i$ which tends to $0$ as $i \rightarrow \infty$
 we have $x_1/N = \sum_{i = 1}^{\infty}a_i/B^i$ which we write as
 $x_1/N = 0.a_1 a_2 \ldots a_i \ldots$.
This is the base $B$ expansion of $x_1/N$ ; $B$ being fixed we omit it from the
 notation.
Reading the equations (\ref{eq2}) mod $N$ shows that for $i \ge 1$
  \begin{equation}\label{eq3}
    x_{i+1} \equiv Bx_i \equiv B^2x_{i-1} \equiv \ldots
    \equiv B^ix_1 \pmod N.
  \end{equation}


Let $e$ be the order of $B \mod N$; denoted $e = \ord(B,N)$.
This means $e$ is the smallest positive integer for which
 $B^e \equiv 1\pmod N$ and $B^f \equiv 1 \pmod N~{\rm iff}~e|f$.
By (\ref{eq3}), $x_{e+1} \equiv B^ex_1 \equiv x_1 \pmod N$ and
 $x_{i+1} \not\equiv x_1 \pmod N$ for $1 \le i < e$.
Since $x_1$, $x_{e + 1}$ both belong to $N^*$, $|x_1 - x_{e + 1}| < N$, so their
 congruence forces $x_{e + 1} = x_1$.
Then $a_{e + 1} = a_1$, $x_{e + 2} = x_2$ and in general $x_{i + e} = x_i$,
 $a_{i + e} = a_i$, $i \ge 1$.
Thus the system (\ref{eq2}) consists of the first $e$ equations which then
 repeat forever.
In particular, the base $B$ expansion of $x_1/N$ is periodic with length $e$
 and we write it as $x_1/N = 0.\overline{a_1 a_2 \ldots a_e}$.
Since $e$ depends only on $N$ and $B$, not $x_1$, we see that every fraction
 $x/N$ with $x \in N^*$ has period length $e$.
Grouping the terms of the infinite series for $x_1/N$ into blocks of $e$ terms
 each, and setting $A = [a_1 a_2 \ldots a_e]_B$, produces the geometric
 series $\sum_{i = 1}^\infty\frac{A}{B^{ei}}$ and shows $x_1/N = A/(B^e - 1)$.
It may be helpful to do a simple numerical example: find the periodic expansion
 of $1/14$ in base 5.
$N = 14$, $B = 5$, $x_1 = 1$; we don't need to know $e = \ord(5,14)$ in
 advance.
The equations (\ref{eq2}) now are
\begin{equation}\label{eq4}
  \begin{matrix}
    5 \cdot  1 = 0 \cdot 14 +  5 \\
    5 \cdot  5 = 1 \cdot 14 + 11 \\
    5 \cdot 11 = 3 \cdot 14 + 13 \\
    5 \cdot 13 = 4 \cdot 14 +  9 \\
    5 \cdot  9 = 3 \cdot 14 +  3 \\
    5 \cdot  3 = 1 \cdot 14 +  1 .
  \end{matrix}
\end{equation}
Having reached the remainder $x_7 = 1 = x_1$, we know that $e = 6$ and
 $1/14 = 0.\overline {013431}$ in base $5$.                 


Let $d$ be a divisor of $e$ and let $k = e/d$, $e = dk$.
Break up the first $e$ equations of (\ref{eq2}) into $d$ groups of $k$
 equations each.
For $1 \le j \le d$, the $j$th group consists of the following $k$ equations.


\begin{equation}\label{eq5}
  \begin{matrix}
    Bx_{(j - 1)k + 1} = a_{(j - 1)k + 1}N + x_{(j - 1)k + 2} \\
    Bx_{(j - 1)k + 2} = a_{(j - 1)k + 2}N + x_{(j - 1)k + 3} \\
    \dots                                                    \\
    Bx_{jk}           = a_{jk}          N + x_{jk + 1} .
  \end{matrix}
\end{equation}


Multiply the first equation by $B^{k-1}$, the second by $B^{k-2}$, $\ldots$,
 the $(k - 1)$th by $B$ and the $k$th by $B^0 = 1$ to obtain


\begin{equation}\label{eq6}
  \begin{matrix}
    B^k    x_{(j-1)k + 1} = a_{(j-1)k + 1}B^{k-1}N + B^{k-1}x_{(j-1)k + 2} \\
    B^{k-1}x_{(j-1)k + 2} = a_{(j-1)k + 2}B^{k-2}N + B^{k-2}x_{(j-1)k + 3} \\
    \dots                                                                  \\
    Bx_{jk}               = a_{jk}               N +        x_{jk + 1} .
  \end{matrix}
\end{equation}


In (\ref{eq6}), the rightmost term of each equation is the left side of the
 next equation; so replace the rightmost term of the first equation by the
 right side of the second equation, then replace the rightmost term of the
 resulting equation by the right side of the third equation, and so on.
Eventually one has
  \begin{equation}\label{eq7}
    B^kx_{(j - 1)k + 1} = (a_{(j - 1)k + 1}B^{k - 1}
    + a_{(j - 1)k + 2}B^{k - 2} + \ldots + a_{jk})N + x_{jk + 1} .
  \end{equation}
The quantity in parentheses is
 $[a_{(j - 1)k + 1} a_{(j - 1)k + 2} \ldots a_{jk}]_B$, the number
 represented by the base $B$ numeral consisting of the $j$th block of $k$
 $B$-digits in the period; denote this number by $A_j$.
So (\ref{eq7}) now becomes $B^kx_{(j - 1)k + 1} = A_jN + x_{jk + 1}$.
Add these equations (\ref{eq7}) for $j = 1$, $2$, $\ldots$, $d$ to obtain
  \begin{equation}\label{eq8}
    B^k\sum_{j = 1}^d x_{(j - 1)k + 1}
    = N(\sum_{j = 1}^d A_j) + \sum_{j = 1}^d x_{jk + 1} .
  \end{equation}
But both sums over $x$ are equal since $x_{dk + 1} = x_{e + 1} = x_1$, so
 (\ref{eq8}) may be rewritten as
  \begin{equation}\label{eq9}
    (B^k - 1)\sum_{j = 1}^d x_{(j - 1)k + 1} = N(\sum_{j = 1}^d A_j) .
  \end{equation}
This relation between the two sums is the key to all that follows.
It is convenient to define
  \begin{equation}\label{eq10}
    R_d(x) = \sum_{j = 1}^d x_{(j - 1)k + 1}
    ~~ {\rm and} ~~ S_d(x) = \sum_{j = 1}^d A_j .
  \end{equation}
Call the set
 $\{x_1,x_{k + 1},\ldots,x_{(d - 1)k + 1}\} = \{x_{jk + 1}~|~j \bmod d\}$ the
 $d$-cycle of $x_1$; more generally, for any $i \ge 1$,
 $\{x_i,x_{k + i},\ldots,x_{(d - 1)k + i}\} = \{x_{jk + i}~|~j \bmod d\}$ is the
 $d$-cycle of $x_i$.
For any two indices $s$ and $t$, $x_s$ and $x_t$ have the same $d$-cycle iff
 $s \equiv t \pmod k$ and for any $x \in N^*$, $R_d(x)$ and $S_d(x)$
 depend only on the $d$-cycle of $x$.   
Of course, $R$ and $S$ depend also on $B$, $N$ and $e = dk$, but we consider
 these fixed for the discussion. 
We summarize the above as 


\begin{theorem}\label{th1}
Given $N$, $B$, $e = \ord(B,N)$ and $e = dk$.
Let $x \in N^*$ and $x/N = 0.\overline{a_1 a_2 \ldots a_e}$ in base $B$.
Break up the period $a_1 a_2 \ldots a_e$ into $d$ blocks of length $k$ each.
For $j = 1$, $2$, $\ldots$, $d$, let $A_j = [a_{(j - 1)k + 1} \ldots a_{jk}]_B$,
 the number represented by the base $B$ numeral consisting of the $j$th  block.
Let $x_1 = x$, $x_2$, $\ldots$ be the remainders in the long division algorithm
 {\rm (\ref{eq2})} for $x/N$.
Then
  \begin{equation}\label{eq11a}
    S_d(x) = (R_d(x)/N)(B^k - 1),
  \end{equation}
  \begin{equation}\label{eq11b}
    S_d(x) \equiv 0 \pmod {B^k - 1} ~~
    iff ~~ R_d(x) \equiv 0 \pmod N.
  \end{equation}
\end{theorem}


\begin{proof}
(\ref{eq11a}) is just a rewriting of (\ref{eq9}) in the notation (\ref{eq10})
 and then (\ref{eq11b}) is immediate.
\end{proof}


\noindent
{\bf Definition.}
Let $N$, $B$, $e$, $d$, $k$ be as above.
We say $N$ has the base $B$ Midy property for the divisor $d$ (of $e$) if
 for every $x \in N^*$, $S_d(x) \equiv 0 \pmod {B^k - 1}$.
We denote by $M_d(B)$ the set of integers that have the Midy property in base
 $B$ for the divisor $d$.


\begin{theorem}\label{th2}
The following are equivalent:
\begin{enumerate}
\item[$(i)$] $N \in M_d(B)$
\item[$(ii)$] For some $x \in N^*$, $S_d(x) \equiv 0 \pmod {B^k - 1}$
\item[$(iii)$] For some $x \in N^*$, $R_d(x) \equiv 0 \pmod N$
\item[$(iv)$] $B^{k(d - 1)} + B^{k(d - 2)} + \ldots + B^k + 1 \equiv 0 \pmod N.$
\end{enumerate}
Furthermore $\gcd(B^k - 1,N) = 1$ implies $N \in M_d(B)$.
\end{theorem}


\begin{proof}
The equivalence of $(ii)$ and $(iii)$ follows from Theorem \ref{th1}.
Noting (\ref{eq3}), we have $$R_d(x) = \sum_{j = 1}^dx_{(j - 1)k + 1}
 \equiv ( \sum_{j = 1}^dB^{(j-1)k})x\pmod N.$$
Since $\gcd(x,N) = 1$,
 $$R_d(x) \equiv 0 \pmod N ~~ {\rm iff} ~~
 \sum_{j = 1}^dB^{k(j - 1)} \equiv 0 \pmod N,$$
 showing $(iv)$ equivalent to $(ii)$ and $(iii)$.
Now $(iv)$ is independent of $x$, so $(iv)$ is equivalent to saying
 $S_d(x) \equiv 0 \pmod {B^k - 1}$ for every $x \in N^*$, which, by
 definition, is $(i)$.
For the last statement, let $F_d(t)$ be the polynomial
 $t^{d - 1} + t^{d - 2} + \ldots + t + 1$, so $(iv)$ amounts to
 $F_d(B^k) \equiv 0 \pmod N$.
But $(B^k - 1)F_d(B^k) = B^e - 1 \equiv 0 \pmod N$, by definition of $e$.
Thus $\gcd(B^k - 1,N) = 1$ implies $(iv)$, hence $N \in M_d(B)$, completing the
 proof.
\end{proof}


Here is an example to show that $\gcd(B^k - 1,N) = 1$ is only sufficient for
 $N \in M_d(B)$, but not necessary.
Take $B = 10$, $N = 21$, $1/21 = 0.\overline{047619}$, $e = 6$.
With $d = 3$, $k = 2$,
 $S_3(1) = 04 + 76 + 19 = 99 \equiv 0 \pmod {10^2 - 1}$, so
 $21 \in M_3(10)$, but $\gcd(10^2 - 1,21) = 3 \ne 1$.


For a numerical illustration, take $N = 14$, $B = 5$, $e = 6$, $x = 1$, as in
 (\ref{eq4}) above.
The period is $013431$ and the remainders $x_1$, $\ldots$, $x_6$ are $1$, $5$,
 $11$, $13$, $9$, $3$, respectively.
With $d = 2$, $k = 3$,
 $S_2(1) = A_1 + A_2 = [013]_5 + [431]_5 = [444]_5 = 5^3 - 1$ and
 $R_2(1) = x_1 + x_4 = 1 + 13 = 14$; thus $14 \in M_2(5)$.
With $d = 3$, $k = 2$,
 $S_3(1) = A_1 + A_2 + A_3 = [01]_5 + [34]_5 + [31]_5 = 36
 \not\equiv 0 \pmod {5^2 - 1}$,
 $R_3(1) = x_1 + x_3 + x_5 = 1 + 11 + 9 = 21$; so $14 \not\in M_3(5)$. 
Note that the relation (\ref{eq11a}) holds: $36 = (21/14)(5^2 - 1)$. 


For $d = 1$, $k = e$, we never have $N \in M_1(B)$, for this would imply
 $1 = R_1(1) \equiv 0 \pmod N$, which is impossible.
Equivalently, $N \in M_1(B)$ says that for any $x \in N^*$,
 $S_1(x) \equiv 0 \pmod {B^e - 1}$.
But $S_1(x) = A = [a_1 a_2 \ldots a_e]_B$, and we've seen that
 $A/(B^e - 1) = x/N$.
So $M_1(B)$ is empty; from now on we consider only $d > 1$.
For $d = e$, $k = 1$, $S_e(x)$ is $\sum_{j = 1}^ea_j$, the sum of the
 $B$-digits in the period.
By Theorem \ref{th2}, $S_e(x) \equiv 0 \pmod {B - 1}$ if $(B - 1,N) = 1$.
In particular, with $B = 10$, the period of the decimal for $x/N$ has the sum
 of its digits divisible by $9$ whenever $N$ is not divisible by $3$.


Given $B$ and $d$ (both $> 1$) it would be nice to be able to describe all
 numbers having the base $B$ Midy property for the divisor $d$; here we make
 only a few observations in this direction.


\begin{theorem}\label{th3}
If $p$ is a prime that does not divide $B$ and $e = \ord(B,p)$ is a multiple of
 $d$, then $p \in M_d(B)$.
Then also $p^h \in M_d(B)$ for every $h > 0$.
\end{theorem}


\begin{proof}
Write $e = dk$; $k < e$ since $d > 1$, so $B^k \not\equiv 1 \pmod p$,
 hence $\gcd(B^k - 1,p) = 1$ and the result follows from Theorem \ref{th2}.
Note that $p$ is not $2$, for if so then $B$ is odd and $B^1 \equiv 1 \pmod 2$,
 so $\ord(B,2) = 1$, which is not a multiple of $d$.
For $p \ne 2$ it is known that $e_h = \ord(B,p^h) = ep^g$, where $g$, depending
 on $h$, is an integer $\ge 0$ whose exact value is not relevant here; see
 \cite[p. 52]{5}.
Thus $e_h = dK$, where $K = kp^g$.
By Fermat, $B^K = (B^k)^{p^g} \equiv B^k \pmod p$, so
 $\gcd(B^K - 1,p^h) = \gcd(B^k - 1,p) = 1$, and the result follows again from
 Theorem \ref{th2}.
\end{proof}


Suppose $p_1$, $p_2$, $\ldots$, $p_r$ are distinct primes all belonging to
 $M_d(B)$ and $N = p_1^{h_1} p_2^{h_2} \ldots p_r^{h_r}$, where $h_1$, $h_2$,
 $\ldots$, $h_r$ are positive integers.
Does $N \in M_d(B)$?
It turns out that the answer does not depend on the values of the the $h_i$.
For $i = 1$, $2$, $\ldots$, $r$, let $\ord(B,p_i) = e_i = dk_i$,
 $E_i = \ord(B,p_i^{h_i}) = e_ip_i^{g_i}$, $g_i \ge 0$.
Now $E = \ord(B,N) = \lcm(E_1,\ldots,E_r)
 = \lcm(dk_1p_1^{g_1},\ldots,dk_rp_r^{g_r})$.
Set $K_i = E_i/d = k_ip_i^{g_i}$ and $K = E/d$, so $K = \lcm(K_1,\ldots,K_r)$.
We need some preliminary remarks.
If $q$ is a prime and $w$ a positive integer, denote by $v_q(w)$ the
 multiplicity of $q$ as a factor of $w$.
Thus
  \begin{equation}\label{eq12}
    w = \prod_qq^{v_q(w)},
    ~~ {\rm the~product~taken~over~all~prime~numbers} ~~ q,
  \end{equation}
where almost all the exponents are $0$.
For positive integers $w_1$, $\ldots$, $w_r$, 
  \begin{equation}\label{eq13}
    \lcm(w_1,\ldots,w_r) = \prod_qq^{m_q},
    ~~ {\rm where} ~~ m_q = \max(v_q(w_1),\ldots,v_q(w_r)).
  \end{equation}
If $Q$ is a set of primes, denote by $Q'$ its complement in the set of all
 primes.
Define the $Q$ part of $w$ to be $u = \prod_{q \in Q}q^{v_q(w)}$ and the $Q'$
 part $y = \prod_{q \in Q'}q^{v_q(w)}$ so by (\ref{eq12}) $w = uy$.
In the same way, (\ref{eq13}) says
 $$\lcm(w_1,\ldots,w_r) = \lcm(u_1,\ldots,u_r) \lcm(y_1,\ldots,y_r).$$
Returning to $N$ above, let $Q$ be the set of primes which divide $d$, and $Q'$
 the complementary set.
Note that each $p_i$ belongs to $Q'$, because $d|e_i \le p_i - 1 < p_i$.
Finally, let $c_i$ be the the largest integer $\ge 0$ for which $d^{c_i}$
 divides $k_i$; so $k_i = d^{c_i}w_i$ and $d \not| ~ w_i$.
Let $u_i$ be the $Q$ part of $w_i$ and $y_i$ the $Q'$ part.
Thus $K_i = k_ip_i^{g_i} = (d^{c_i}u_i)(y_ip_i^{g_i})$ is the factorization of
 $K_i$ into the product of its $Q$ part and $Q'$ part, and
 $$K = \lcm(K_1,\ldots,K_r)
 = \lcm(d^{c_1}u_1,\ldots,d^{c_r}u_r) \lcm(y_1p_1^{g_1},\ldots,y_rp_r^{g_r}).$$
Set
  \begin{equation}\label{eq14}
    U = \lcm(d^{c_1}u_1,\ldots,d^{c_r}u_r),
    ~~ Y = \lcm(y_1p_1^{g_1},\ldots,y_rp_r^{g_r}).
  \end{equation}
so $K = UY$ is the factorization of $K$ into the product of its $Q$ part $U$
 and $Q'$ part $Y$.


\begin{theorem}\label{th4}
Let $p_1$, $\ldots$, $p_r$ be primes each belonging to $M_d(B)$ and $h_1$,
 $\ldots$, $h_r$ positive integers and $N = p_1^{h_1} \ldots p_r^{h_r}$.
With the notations introduced above, $N \in M_d(B)$ if and only if 
  \begin{equation}\label{eq15}
    for ~~ i = 1,\ldots,r, ~~~ U/(d^{c_i}u_i) \not\equiv 0 \pmod d.                                                      
  \end{equation}
This condition depends only on the primes $p_1$, $\ldots$, $p_r$ and not the
 exponents $h_1$, $\ldots$, $h_r$.
If $d$ is a prime $q$, $N \in M_q(B)$ if and only if $q$ occurs with the same
 multiplicity in each $e_i$:
  \begin{equation}\label{eq16}
    v_q(e_1) = v_q(e_2) = \ldots = v_q(e_r).
  \end{equation}
\end{theorem}


\begin{proof}
Clearly, by definition of $U$, $U/(d^{c_i}u_i)$ is an integer for each $i$.
If for some $i$, $U/(d^{c_i}u_i) \equiv 0 \pmod d$ then $dd^{c_i}u_i|U$
 and, since also $y_ip_i^{g_i}|Y$, it follows that
 $E_i = dd^{c_i}u_iy_ip_i^{g_i}|UY = K$.
Hence $B^K \equiv 1 \pmod {p_i^{h_i}}$ and, in particular,
 $B^K \equiv 1 \pmod {p_i}$.
Then $F_d(B^K) = \sum_{j = 1}^d(B^K)^{j - 1} \equiv \sum_{j = 1}^d1 \equiv d
 \pmod {p_i}$.
Now by Theorem \ref{th2}, if $N \in M_d(B)$ then $F_d(B^K) \equiv 0 \pmod N$
 implying $F_d(B^K) \equiv 0 \pmod {p_i}$, which combined with the previous
 congruence shows $d \equiv 0 \pmod {p_i}$ which is absurd since $d|e_i < p_i$.
So the condition (\ref{eq15}) is necessary for $N \in M_d(B)$.
Suppose now that (\ref{eq15}) is satisfied.
Then for each $i$, $dd^{c_i}u_i \not| ~ U$, so
 $e_i = dd^{c_i}u_iy_i \not| ~ UY = K$ and so $B^K \not\equiv 1 \pmod {p_i}$.
Thus for each $i$, $(B^K - 1,p_i) = 1$, hence $(B^K - 1,N) =1$, which, by
 Theorem \ref{th2}, implies $N \in M_d(B)$.
This proves (\ref{eq15}) is also sufficient, and clearly (\ref{eq15}) is
 independent of $h_1$, $\ldots$, $h_r$.
Now consider the case where $d$ is a prime number $q$, then $Q = \{q\}$
 consists of the single prime $q$.
Then the definition of $c_i$ as the the largest integer for which $q^{c_i}|k_i$
 says $c_i = v_q(k_i)$; thus $k_i = q^{c_i}w_i$ and $q \not| ~ w_i$ so the $Q$
 part $u_i$ of $w_i$ is $1$ which means
 $U = lcm(q^{c_1},\ldots,q^{c_r}) = q^c$, where $c = max(c_1,\ldots,c_r)$.
Hence the conditions of (\ref{eq15}) become simply that for each $i$,
 $q^c/q^{c_i}$ is not divisible by $q$, so $c_i = c$ and
 $v_q(e_i) = v_q(qk_i) = 1 + c$.
This completes the proof.
\end{proof}


Theorem \ref{th4} was first proved by Jenkins \cite{3} in the case $d = 2$.
In \cite[p. 94]{7}, the author seems to claim that if $d$ is any integer, prime
 or not, then $N \in M_d(B)$ if and only if all the $c_i$ are equal:
 $c_1 = c_2 = \ldots = c_r$.
As our proof shows, this is true when $d$ is a prime but not otherwise.


Here are numerical illustrations of some of our results, which will also show
 that the above claim is false.
We keep the usual notations.
Let $p_1 = 7$, $p_2 = 9901$, $p_3 = 19$, $B = 10$: $1/7 = 0.\overline{142857}$,
 $e_1 = 6$; $1/9901 = 0.\overline{000100999899}$, $e_2 = 12$;
 $1/19 = 0.\overline{052631578947368421}$, $e_3 = 18$.
One checks easily that each $p_i \in M_d(10)$, for each $d | 6$, $d > 1$, as
 stated in Theorem \ref{th3}.
For example, for $19$ with $d = 6$, $k = 3$,
 $S_6(1) = 052 + 631 + 578 + 947 + 368 + 421 = 2997 \equiv 0 \pmod {10^3 - 1}$.
Note that in the setup of Theorem \ref{th4}, whenever some $h_i = 1$, then
 $g_i = 0$, $E_i = e_i$, $K_i = k_i$; this will be the case in what follows.
Now for $p_1p_2 = 7 \times 9901 = 69307$, $E = \lcm(6,12) = 12$,
 $1/69307 = 0.\overline{000014428557}$.
Consider, for Theorem \ref{th4}, those $d$ which divide both $6$ and $12$: $2$,
 $3$, which are primes, and $6$ which is not.
$v_3(6) = 1 = v_3(12)$, so $69307 \in M_3(10)$, while
 $v_2(6) = 1 \ne v_2(12) = 2$, so $69307 \not\in M_2(10)$, as one also easily
 verifies from the period.
For $d = 6$, $Q = \{2,3\}$, $K_1 = k_1 = 1$, $K_2 = k_2 = 2$, $c_1 = c_2 = 0$,
 $u_1 = y_1 = 1$, $u_2 = 2$, $y_2 = 1$, and (\ref{eq14}) gives
 $U = \lcm(1,2) = 2$, $Y = \lcm(1,1) = 1$, $K = UY = 2$.
(\ref{eq15}) is satisfied: for $i = 1$, $2/1 \not\equiv 0 \pmod 6$; for
 $i = 2$, $2/2 \not\equiv 0 \pmod 6$.
Thus we know $69307 \in M_6(10)$; again we verify this directly from the
 period.
$S_6(1) = 00 + 00 + 14 + 42 + 85 + 57 = 198 \equiv 0 \pmod {10^2 - 1}$.
For a later application we note here that
 $S_4(1) = 000 + 014 + 428 + 557 = 999$.


Now consider $p_1p_2p_3 = 7 \times 9901 \times 19 = 1316833$,
 $E = \lcm(6,12,18) = 36$,
 $1/1316833 = 0.\overline{000000759397736842864660894737601503}$.
For $d = 2$, $v_2(6) = 1$, $v_2(12) = 2$, $v_2(18) = 1$ and for $d = 3$,
 $v_3(6) = v_3(12) = 1$, $v_3(18) = 2$, so $1316833$ is not in $M_d(10)$ for
 $d = 2$ and $3$--again this can be verified from the period.
With $d = 6$, $K_1 = k_1 = 1$, $K_2 = k_2 = 2$, $K_3 = k_3 = 3$; none of these
 is divisible by $6$, so $c_1 = c_2 = c_3 = 0$.
$Q = \{2,3\}$, $u_1 = 1$, $u_2 = 2$, $u_3 = 3$ while $y_1 = y_2 = y_3 = 1$,
 $U = \lcm(1,2,3) = 6$, $Y = 1$, $K = 6$.
Now consider (\ref{eq15}): for $i = 1$,
 $U/(d^{c_1}u_1) = 6/1 \equiv 0 \pmod 6$, so the condition is not satisfied and
 $1316833 \not\in M_6(10)$.
This is a counterexample to the aforementioned claim.
To check this numerically, $S_6(1)
 = 000000 + 759397 + 736842 + 864660 + 894737 + 601503
 = 3857139 \not\equiv 0 \pmod {10^6 - 1}$.
In fact, $3857139/999999 = 27/7$.
The other divisors of $36$ which do not arise from Theorem \ref{th4} are
 $d = 4$, $9$, $12$, $18$, $36$ and the reader may verify that
 $1316833 \in M_d(10)$ for each of these.
The next theorem shows that not all of this is accidental, but that once it is
 known for $4$ and $9$ the result follows for their multiples $12$, $18$, $36$.


\begin{theorem}\label{th5}
Suppose $e = \ord(B,N)$ and $d_1|d_2$, $d_2|e$.
If $N \in M_{d_1}(B)$ then $N \in M_{d_2}(B)$.
\end{theorem}


\begin{proof}
Write $e = d_1k_1 = d_2k_2$ and set $c = d_2/d_1 = k_1/k_2$.
Since $N \in M_{d_1}(B)$, $R_{d_1}(x) \equiv 0 \pmod N$ for every $x \in N^*$.
By definition, $R_{d_2}(x) = \sum_{j = 0}^{d_2 - 1}x_{jk_2 + 1}$.
We will show that $R_{d_2}(x) = \sum_{r = 0}^{c - 1}R_{d_1}(x_{rk_2 + 1})$,
 hence $R_{d_2}(x)$ is a sum of terms $\equiv 0 \pmod N$ so it is also
 $\equiv 0 \pmod N$ which implies $S_{d_2}(x) \equiv 0 \pmod {B^{k_2} - 1}$
 and $N \in M_{d_2}(B)$.
The numbers $j = 0$, $1$, $\dots$, $d_2 - 1 = cd_1 - 1$ may be written as
 $j = ic + r$, where $i = 0$, $1$, $\ldots$, $d_1 - 1$ and $r = 0$, $1$,
 $\ldots$, $c - 1$; then $jk_2 + 1 = ick_2 + rk_2 + 1 = ik_1 + rk_2 + 1$.
Thus $$R_{d_2}(x)
 = \sum_{r = 0}^{c - 1}\sum_{i = 0}^{d_1 - 1}x_{ik_1 + rk_2 + 1},$$
 and the inner sum is just $R_{d_1}(x_{rk_2 + 1})$; this completes the proof.
\end{proof}


The basic idea here is that the $d_2$-cycle of $x$ is a union of $c$
 $d_1$-cycles.


\section{The Multiplier}\label{s3}


For $N \in M_d(B)$ and $x \in N^*$ we have, by definition,
 $S_d(x) \equiv 0 \pmod {B^k - 1}$, and more precisely, by (\ref{eq11a}),
 $S_d(x) = m_d(x)(B^k - 1)$ where $m_d(x) = R_d(x)/N$ is an integer, which we
 call the multiplier; in general it depends on both $d$ and the $d$-cycle of
 $x$.


\begin{theorem}\label{th6}
If $N \in M_2(B)$ then for every even $d|e$, $N \in M_d(B)$ and
 $m_d(x) = d/2$ for every $x \in N^*$.
\end{theorem}


\noindent
\begin{remark}
Midy's Theorem of the {\rm Introduction} now follows.
For taking $B = 10$, the conditions stated there about $N$ show, by Theorems
 \ref{th3} and \ref{th2}, that $N \in M_2(10)$ and then this Theorem shows
 $m_2(x) = 1$, so $S_2(x) = 10^k - 1$, which is a string of $k = e/2$ 9's.
\end{remark}


\begin{proof}
Let $e = 2k$.
By Theorem \ref{th2}$(iv)$, $B^k + 1 \equiv 0$, or $B^k \equiv -1 \pmod N$,
 which, by (\ref{eq3}) with $i = k + 1$, shows $x_{k + 1} \equiv -x_1 \pmod N$.
But the only member of $N^*$ that is congruent to $-x_1$ is $N - x_1$, hence
 $x_{k+1} = N - x_1$, so $R_2(x) = x_1 + x_{k + 1} = x_1 + (N - x_1) = N$,
 $m_2(x) = 1$; this proves the case $d = 2$.
Now say $d > 2$, $2 | d$, $d | e$, $c = d/2$, $k' = e/d$; as shown in the
 proof of Theorem \ref{th5} the $d$-cycle of $x$ is a union of $c$ $2$-cycles
 and $R_d(x) = \sum_{r = 0}^{c - 1}R_2(x_{rk' + 1}) = \sum_{r = 0}^{c - 1}N
 = (d/2)N$, hence $m_d(x) = d/2$.
\end{proof}


The condition $N \in M_2(B)$ in Theorem \ref{th6} cannot be omitted.
For example, we've seen--after the proof of Theorem \ref{th4}--that for
 $N = 69307$, $e = 12$, $N$ does not belong to $M_2(10)$ but $N$ does belong to
 $M_4(10)$ and $S_4(1) = 999 = (10^3 - 1)$.
Thus in this case $d = 4$ is even and $m_4(1) = 1 \ne 4/2$.


We now study the multiplier $m_3(x)$  for $N \in M_3(B)$, $e = 3k$.
Recall the result of Ginsberg \cite{2} stated in the Introduction which, in our
 current notation, says $m_3(1) = 1$ if $N$ is a prime.
We now show that such a result holds much more extensively.


\begin{theorem}\label{th7}
Suppose $N \in M_3(B)$, $e = 3k$. Then
\begin{enumerate}
\item[$(i)$]   $m_3(1) = 1$
\item[$(ii)$]  if $N$ is odd, $m_3(2) = 1$
\item[$(iii)$] if $3 \not| ~ N$ and $N \ne 7$, $m_3(3) = 1$.
\end{enumerate}
\end{theorem}


\begin{proof}
For $x \in N^*$, $R_3(x) = x_1 + x_{k + 1} + x_{2k + 1} < N + N + N = 3N$.
Since $R_3(x) \equiv 0 \pmod N$, $R_3(x_1) = N$ or $2N$.
If $x = 1$ or $2$, then $x_{k + 1}$, $x_{2k + 1}$ are at most $N - 1$,
 $N - 2$ (in some order).
Thus $R_3(x) \le 2 + (N - 1) + (N - 2) < 2N$, which forces $R_3(x) = N$,
 $m_3(x_1) = 1$, proving $(i)$ and $(ii)$.
Now take $x = 3$; $R_3(3) \le 3 + (N - 1) + (N - 2) \le 2N$, where equality
 holds iff $x_{k + 1} = N - 1$, $x_{2k + 1} = N - 2$, or $x_{k + 1} = N - 2$,
 $x_{2k + 1} = N - 1$.
In the former case, by (\ref{eq3}), $N - 1 \equiv 3B^k$ and
 $N - 2 \equiv 3B^{2k} \pmod N$, so $9B^{3k} \equiv 2 \pmod N$.
But $3k = e$, $B^{3k} \equiv 1 \pmod N$, so $9 \equiv 2 \pmod N$, hence
 $N = 7$.
In the latter case the argument is the same with $N - 1$, $N - 2$ interchanged.
This proves $(iii)$.
\end{proof}


Note that $7$  really is exceptional; take, say, $B = 10$,
 $3/7 = 0.\overline{428571}$, $S_3(3) = 42 + 85 + 71 = 198 = 2(10^2 - 1)$, so
 here $m_3(3) = 2$.


\section{Conclusion}\label{s4}


Midy's Theorem and its extensions deserve to be better known and certainly
 have a place in elementary number theory.
These patterns in the decimal expansions of rational numbers provide an
 unexpected glimpse of the charm, and structure, of mathematical objects.
Many questions and unexplored pathways remain to be investigated.



\noindent
Department of Mathematics and Computer Science,
 Lehman College, CUNY, Bronx, NY 10468


\noindent
joseph.lewittes@lehman.cuny.edu


\end{document}